\input amstex
\input amsppt.sty
\magnification=\magstep1
\hsize=30truecc
\vsize=22.2truecm
\baselineskip=16truept
\NoBlackBoxes
\TagsOnRight \pageno=1 \nologo
\def\Z{\Bbb Z}
\def\N{\Bbb N}

\def\l{\left}
\def\r{\right}
\def\bg{\bigg}
\def\({\bg(}
\def\[{\bg\lfloor}
\def\){\bg)}
\def\]{\bg\rfloor}
\def\t{\text}
\def\f{\frac}

\def\bi{\binom}
\def\eq{\equiv}

\def\ls{\leqslant}
\def\gs{\geqslant}
\def\mo{\roman{mod}}

\def\da{\delta}

\def\Proof{\noindent{\it Proof}}

\def\Remark{\medskip\noindent{\it  Remark}}

\hbox {J. Number Theory 131(2011), no.\,12, 2387--2397.}
\bigskip
\topmatter
\title On Delannoy numbers and Schr\"oder numbers\endtitle
\author Zhi-Wei Sun\endauthor
\leftheadtext{Zhi-Wei Sun}
\rightheadtext{On Delannoy numbers and Schr\"oder numbers}
\affil Department of Mathematics, Nanjing University\\
 Nanjing 210093, People's Republic of China
  \\  zwsun\@nju.edu.cn
  \\ {\tt http://math.nju.edu.cn/$\sim$zwsun}
\endaffil
\abstract The $n$th Delannoy number and the $n$th Schr\"oder number  given by
$$D_n=\sum_{k=0}^n\bi nk\bi{n+k}k\ \ \t{and}\ \ S_n=\sum_{k=0}^n\bi nk\bi{n+k}k\f1{k+1}$$
respectively arise naturally from enumerative combinatorics. Let $p$ be an odd prime. We mainly show that
$$\sum_{k=1}^{p-1}\f{D_k}{k^2}\eq2\l(\f{-1}p\r)E_{p-3}\ (\mo\ p)$$
and
$$\sum_{k=1}^{p-1}\f{S_k}{m^k}\eq\f{m^2-6m+1}{2m}\(1-\l(\f{m^2-6m+1}p\r)\)\ (\mo\ p),$$
where $(-)$ is the Legendre symbol, $E_0,E_1,E_2,\ldots$ are Euler numbers, and $m$ is any integer not divisible by $p$.
We also conjecture that
$$\sum_{k=1}^{p-1}\f{D_k^2}{k^2}\eq -2q_p(2)^2\ (\mo\ p)$$
 where $q_p(2)$ denotes the Fermat quotient $(2^{p-1}-1)/p$.
\endabstract
\thanks 2010 {\it Mathematics Subject Classification}.\,Primary 11A07, 11B75;
Secondary  05A15, 11B39, 11B68, 11E25.
\newline\indent Supported by the National Natural Science
Foundation (grant 10871087) of China.
\endthanks
\endtopmatter
\document

\heading{1. Introduction}\endheading

For $n\in\N=\{0,1,2,\ldots\}$, the (central) Delannoy number $D_n$ denotes the
number of lattice paths from the point $(0,0)$ to $(n,n)$ with steps
$(1,0),(0,1)$ and $(1,1)$, while the Schr\"oder number $S_n$
represents the number of such paths that never rise above the line
$y=x$. It is known that
$$D_n=\sum_{k=0}^n\bi nk\bi{n+k}k=\sum_{k=0}^n\bi {n+k}{2k}\bi{2k}k$$
and
$$S_n=\sum_{k=0}^n\bi nk\bi{n+k}k\f1{k+1}=\sum_{k=0}^n\bi {n+k}{2k}C_k,$$
where $C_k$ stands for the Catalan number $\bi{2k}k/(k+1)=\bi{2k}k-\bi{2k}{k+1}$.
For information on $D_n$ and $S_n$, the reader may consult [CHV], [S], and  p.\,178 and p.\,185 of [St].

Despite their combinatorial backgrounds, surprisingly Delannoy numbers and Schr\"oder
numbers have some nice number-theoretic properties.

As usual, for an odd prime $p$ we let $(\f{\cdot}p)$ denote the
Legendre symbol. Recall that Euler numbers $E_0,E_1,E_2,\ldots$ are
integers defined by $E_0=1$ and the recursion:
$$ \sum^n\Sb k=0\\2\mid k\endSb \bi nk E_{n-k}=0\ \ \ \t{for}\ n=1,2,3,\ldots.$$

Our first theorem is concerned with Delannoy numbers and their generalization.

\proclaim{Theorem 1.1} Let $p$ be an odd prime. Then
$$\sum_{k=1}^{p-1}\f{D_k}{k^2}\eq2\l(\f{-1}p\r)E_{p-3}\ (\mo\ p)\tag1.1$$
and
$$\sum_{k=1}^{p-1}\f{D_k}k\eq-q_p(2)\ (\mo\ p),\tag1.2$$
where $q_p(2)$ denotes the Fermat quotient $(2^{p-1}-1)/p$.
If we set
$$D_n(x)=\sum_{k=0}^n\bi nk\bi{n+k}kx^k\ \ (n\in\N),$$
then for any $p$-adic integer $x$ we have
$$\sum_{k=1}^{p-1}\f{D_k(x)}k\eq\f{(-1+\sqrt{-x})^p+(-1-\sqrt{-x})^p+2}p\ (\mo\ p).\tag1.3$$
\endproclaim
\proclaim{Corollary 1.1} Let $p$ be an odd prime. We have
$$\align\sum_{k=1}^{p-1}\f{D_k(3)}k\eq&-2q_p(2)\ \ (\mo\ p)\ \ \t{provided}\ p\not=3,\tag1.4
\\\sum_{k=1}^{p-1}\f{D_k(-4)}k\eq&\f{3-3^p}p\ (\mo\ p),\tag1.5
\\\sum_{k=1}^{p-1}\f{D_k(-9)}k\eq&-6q_p(2)\ (\mo\ p),\tag1.6
\endalign$$
and also
$$\sum_{k=1}^{p-1}\f{D_k(-2)}k\eq-\f 4p P_{p-(\f2p)}\ (\mo\ p),\tag1.7$$
where the Pell sequence $\{P_n\}_{n\gs0}$ is given by
$$P_0=0,\ P_1=1,\ \t{and}\ P_{n+1}=2P_n+P_{n-1}\ (n=1,2,3,\ldots).$$
If $p\not=5$, then
$$\sum_{k=1}^{p-1}\f{D_k(-5)}k\eq-2q_p(2)-\f 5pF_{p-(\f p5)}\ (\mo\ p),\tag1.8$$
where the Fibonacci sequence $\{F_n\}_{n\gs0}$ is defined by
$$F_0=0,\ F_1=1,\ \t{and}\ F_{n+1}=F_n+F_{n-1}\ (n=1,2,3,\ldots).$$
\endproclaim

Now we propose two conjectures which seem challenging in the author's opinion.

\proclaim{Conjecture 1.1} Let $p>3$ be a prime. We have
$$\align\sum_{k=1}^{p-1}\f{D_k^2
}{k^2}\eq&-2q_p(2)^2\ (\mo\ p),\tag1.9
\\\sum_{k=1}^{p-1}\f{D_k}k\eq&-q_p(2)+p\,q_p(2)^2\ (\mo\ p^2),\tag1.10
\\\sum_{k=1}^{p-1}D_kS_k\eq&-2p\sum_{k=1}^{p-1}\f{(-1)^k+3}k\ (\mo\ p^4),
\endalign$$ and
$$\sum_{k=1}^{(p-1)/2}D_kS_k\eq\cases 4x^2\ (\mo\ p)&\t{if}\ p\eq1\ (\mo\ 4)\ \&\ p=x^2+y^2\ (2\nmid x,\, 2\mid y),
\\0\ (\mo\ p)&\t{if}\ p\eq3\ (\mo\ 4).\endcases$$
Also, $\sum_{n=1}^{p-1}s_n^2/n\eq-6\ (\mo\ p)$, where
$$s_n:=\sum_{k=0}^n\bi{n+k}{2k}\bi{2k}{k+1}=D_n-S_n.$$
\endproclaim
\noindent{\it Remark} 1.1. Let $p$ be an odd prime. Though there are
many congruences for $q_p(2)$ mod $p$, (1.9) is curious since its
left-hand side is a sum of squares.  It is known that
$\sum_{k=1}^{p-1}1/k\eq-p^2B_{p-3}/3\ (\mo\ p^3)$ if $p>3$, where
$B_0,B_1,B_2,\ldots$ are Bernoulli numbers. In addition, we can
prove that $\sum_{k=0}^{p-1}D_k\eq(\f{-1}p)-p^2E_{p-3}\ (\mo\ p^3)$
and $\sum_{k=0}^{p-1}D_k^2\eq(\f 2p)\ (\mo\ p)$.

\proclaim{Conjecture 1.2} Let $p>3$ be a prime. Then
$$\align&\sum_{k=0}^{p-1}(-1)^kD_k(2)^3\eq\sum_{k=0}^{p-1}(-1)^k D_k\l(-\f14\r)^3
\eq\l(\f{-2}p\r)\sum_{k=0}^{p-1}(-1)^kD_k\l(\f18\r)^3
\\\eq&\cases(\f{-1}p)(4x^2-2p)\ (\mo\ p^2)&\t{if}\ p\eq1\ (\mo\ 3)\ \&\ p=x^2+3y^2\, (x,y\in\Z),
\\0\ (\mo\ p^2)&\t{if}\ p\eq2\ (\mo\ 3).\endcases
\endalign$$
Also,
$$\align&\l(\f{-1}p\r)\sum_{k=0}^{p-1}(-1)^kD_k\l(\f12\r)^3
\\\eq&\cases4x^2-2p\ (\mo\ p^2)&\t{if}\ p\eq1,7\ (\mo\ 24)\ \t{and}\ p=x^2+6y^2\,(x,y\in\Z),
\\8x^2-2p\ (\mo\ p^2)&\t{if}\ p\eq5,11\ (\mo\ 24)\ \t{and}\ p=2x^2+3y^2\,(x,y\in\Z),
\\0\ (\mo\ p^2)&\t{if}\ (\f {-6}p)=-1.\endcases
\endalign$$
And
$$\align&\sum_{k=0}^{p-1}(-1)^kD_k(-4)^3\eq\l(\f{-5}p\r)\sum_{k=0}^{p-1}(-1)^kD_k\l(-\f1{16}\r)^3
\\\eq&\cases4x^2-2p\ (\mo\ p^2)&\t{if}\ p\eq1,4\ (\mo\ 15)\ \t{and}\ p=x^2+15y^2\,(x,y\in\Z),
\\12x^2-2p\ (\mo\ p^2)&\t{if}\ p\eq2,8\ (\mo\ 15)\ \t{and}\ p=3x^2+5y^2\,(x,y\in\Z),
\\0\ (\mo\ p^2)&\t{if}\ (\f{-15}p)=-1.\endcases
\endalign$$
\endproclaim
\noindent{\it Remark} 1.2. Note that $(-1)^nD_n(x)=D_n(-x-1)$ for
any $n\in\N$, since
$$\align D_n(-x-1)=&\sum_{k=0}^n\bi nk\bi{-n-1}k\sum_{j=0}^k\bi kjx^j
\\=&\sum_{j=0}^n\binom njx^j\sum_{k=0}^n\bi{-n-1}k\binom{n-j}{n-k}
\\=&\sum_{j=0}^n\binom njx^j\binom{-j-1}n=(-1)^nD_n(x).
\endalign$$
\medskip

Concerning Schr\"oder numbers we establish the following result.

\proclaim{Theorem 1.2} Let $p$ be an odd prime and let $m$ be an integer not divisible by $p$. Then
$$\sum_{k=1}^{p-1}\f{S_k}{m^k}\eq\f{m^2-6m+1}{2m}\(1-\l(\f{m^2-6m+1}p\r)\)\ (\mo\ p).\tag1.11$$
\endproclaim
\medskip

\noindent{\it Example}\ 1.1. Theorem 1.2 in the case $m=6$ gives that
$$\sum_{k=1}^{p-1}\f{S_k}{6^k}\eq0\ (\mo\ p)\qquad \ \t{for any prime}\ p>3.\tag1.12$$
\medskip

For technical reasons, we will prove Theorem 1.2 in the next section
and show Theorem 1.1 and Corollary 1.1 in Section 3.

\heading{2. Proof of Theorem 1.2}\endheading

\proclaim{Lemma 2.1} Let $p$ be an odd prime and let $m$ be any
integer not divisible by $p$. Then
$$\sum_{k=1}^{p-1}\f{C_k}{m^k}\eq\f{m-4}2\(1-\l(\f{m(m-4)}p\r)\)\
(\mo\ p).\tag2.1$$
\endproclaim
\Proof. This follows from [Su10, Theorem 1.1] in which the author even determined
$\sum_{k=1}^{p-1}C_k/m^k$ mod $p^2$. However, we will give here a simple proof of (2.1).

For each $k=1,\ldots,p-1$, we clearly have
$$\bi{(p-1)/2}k\eq\bi{-1/2}k=\f{\bi{2k}k}{(-4)^k}\ (\mo\ p).$$
Note also that
$$C_{p-1}=\f1{2p-1}\prod_{k=1}^{p-1}\f{p+k}k\eq-1\ (\mo\ p).$$
Therefore
$$\align\sum_{k=1}^{p-1}\f{C_k}{m^k}\eq&\sum_{0<k<p-1}\bi{(p-1)/2}k\f1{k+1}\l(-\f 4m\r)^k+\f{C_{p-1}}{m^{p-1}}
\\\eq&-\f m4\times\f2{p+1}\sum_{k=1}^{(p-1)/2}\bi{(p+1)/2}{k+1}\l(-\f 4m\r)^{k+1}-1
\\\eq&-\f m2\(\l(1-\f 4m\r)^{(p+1)/2}-1-\f{p+1}2\l(-\f 4m\r)\)-1
\\\eq&-\f m2\(\f{m-4}m\times\f{(m(m-4))^{(p-1)/2}}{m^{p-1}}-1+\f2m\)-1
\\\eq&-\f{m-4}2\l(\f{m(m-4)}p\r)+\f m2-2\ (\mo\ p)
\endalign$$
and hence (2.1) follows.
\qed

\proclaim{Lemma 2.2} For any odd prime $p$ we have
$$\sum_{k=1}^{p-1}S_k\eq2\l(\f{-1}p\r)-2^p\ (\mo\ p^2).\tag2.2$$
\endproclaim
\Proof. Recall the known identity (cf. (5.26) of [GKP, p.\,169])
$$\sum_{n=0}^m\bi nk=\bi{m+1}{k+1}\ \ \ (k,m\in\N).$$
Then
$$\align\sum_{n=0}^{p-1}S_n=&\sum_{n=0}^{p-1}\sum_{k=0}^n\bi{n+k}{2k}C_k=\sum_{k=0}^{p-1}C_k\sum_{n=k}^{p-1}\bi{n+k}{2k}
\\=&\sum_{k=0}^{p-1}C_k\bi{p+k}{2k+1}=\sum_{k=0}^{p-1}\f p{k!(k+1)!(2k+1)}\prod_{0<j\ls k}(p^2-j^2)
\\\eq&\sum_{k=0}^{p-1}\f {p(-1)^k(k!)^2}{k!(k+1)!(2k+1)}
=p\sum_{k=0}^{p-1}(-1)^k\l(\f2{2k+1}-\f1{k+1}\r)
\ (\mo\ p^2).
\endalign$$

Observe that
$$\align 2p\sum_{k=0}^{p-1}\f{(-1)^k}{2k+1}=&p\sum_{k=0}^{p-1}\(\f{(-1)^k}{2k+1}+\f{(-1)^{p-1-k}}{2(p-1-k)+1}\)
\\=&p\sum_{k=0}^{p-1}(-1)^k\l(\f1{2k+1}+\f1{2p-(2k+1)}\r)
\\\eq&p(-1)^{(p-1)/2}\l(\f1p+\f1{2p-p}\r)=2\l(\f{-1}p\r)\ (\mo\ p^2).
\endalign$$
Also,
$$\align &-p\sum_{k=0}^{p-1}\f{(-1)^k}{k+1}=p\sum_{k=1}^p\f{(-1)^k}k
\\\eq&-\sum_{k=1}^{p-1}\f pk\bi{p-1}{k-1}-1=-\sum_{k=0}^{p-1}\bi pk=1-2^p\ (\mo\ p^2).
\endalign$$
Combining the above, we obtain
$$\sum_{n=0}^{p-1}S_n\eq2\l(\f{-1}p\r)+1-2^p\ (\mo\ p^2)$$
and hence (2.2) holds. \qed

\medskip
\noindent{\it Proof of Theorem 1.2}. In the case $m\eq1\ (\mo\ p)$, (1.11) reduces to the congruence
$$\sum_{k=1}^{p-1}S_k\eq-2\l(1-\l(\f{-1}p\r)\r)\ (\mo\ p)$$
which follows from (2.2) in view of Fermat's little theorem.

Below we assume that $m\not\eq1\ (\mo\ p)$. Then
$$\sum_{n=1}^{p-1}\f1{m^n}\eq\sum_{n=1}^{p-1}m^{p-1-n}=\f{m^{p-1}-1}{m-1}\eq0\ (\mo\ p)$$
and hence
$$\sum_{n=1}^{p-1}\f{S_n}{m^n}\eq\sum_{n=1}^{p-1}\f{S_n-1}{m^n}=\sum_{n=1}^{p-1}\f{\sum_{k=1}^n\bi{n+k}{2k}C_k}{m^n}
=\sum_{k=1}^{p-1}\f{C_k}{m^k}\sum_{n=k}^{p-1}\f{\bi{n+k}{2k}}{m^{n-k}}\ (\mo\ p).$$
Given $k\in\{1,\ldots,p-1\}$, we have
$$\sum_{n=k}^{p-1}\f{\bi{n+k}{2k}}{m^{n-k}}=\sum_{r=0}^{p-1-k}\f{\bi{2k+r}r}{m^r}
=\sum_{r=0}^{p-1-k}\f{\bi{-2k-1}r}{(-m)^r}\eq\sum_{r=0}^{p-1-k}\f{\bi{p-1-2k}r}{(-m)^r}\ (\mo\ p).$$
If $(p-1)/2<k<p-1$, then
$$C_k=\f{(2k)!}{k!(k+1)!}\eq0\ (\mo\ p).$$
Therefore
$$\align \sum_{n=1}^{p-1}\f{S_n}{m^n}\eq&\sum_{k=1}^{(p-1)/2}\f{C_k}{m^k}\l(1-\f1m\r)^{p-1-2k}+\f{C_{p-1}}{m^{p-1}}
\\\eq&\sum_{k=1}^{p-1}\f{C_k}{m^k}\l(\f m{m-1}\r)^{2k}\eq\sum_{k=1}^{p-1}\f{C_k}{m_0^k}\ (\mo\ p),
\endalign$$
where $m_0$ is an integer with $m_0\eq(m-1)^2/m\ (\mo\ p)$.
By Lemma 2.1,
$$\align\sum_{k=1}^{p-1}\f{C_k}{m_0^k}\eq&\f{m_0-4}2\(1-\l(\f{m_0(m_0-4)}p\r)\)
\\=&\f{mm_0-4m}{2m}\(1-\l(\f{mm_0(mm_0-4m)}p\r)\)
\\\eq&\f{(m-1)^2-4m}{2m}\(1-\l(\f{(m-1)^2-4m}p\r)\)\ (\mo\ p).
\endalign$$
So (1.11) follows. We are done. \qed

\heading{3. Proofs of Theorem 1.1 and Corollary 1.1}\endheading
\medskip

We need some combinatorial identities.

\proclaim{Lemma 3.1} For any $n\in\N$, we have
$$\sum_{r=0}^{2n}\f{(-1)^r\bi{2n}r}{2n+1-2r}=\f{(-16)^n}{(2n+1)\bi{2n}n}\tag3.1$$
and
$$\sum_{r=0}^{2n}\f{(-1)^r\bi{2n}r}{(2n+1-2r)^2}=\f{(-16)^n}{(2n+1)^2\bi{2n}n},\tag3.2$$
that is,
$$\sum_{k=-n}^n\f{(-1)^k}{(2k+1)^s}\bi{2n}{n-k}=\f{16^n}{(2n+1)^s\bi{2n}n}\ \ \ \t{for}\ s=1,2.\tag3.3$$
\endproclaim
\Proof. If we denote by $a_n$ the left-hand side of (3.1), then the well-known Zeilberger
algorithm (cf. [PWZ]) yields the recursion
$$a_{n+1}=-\f{8(n+1)}{2n+3} a_n\ \ (n=0,1,2,\ldots).$$
So (3.1) can be easily proved by induction.
(3.2) is equivalent to [Su11, (2.5)] which was shown by a similar method.
Clearly (3.3) is just a combination of (3.1) and (3.2). We are done. \qed

\medskip \noindent{\it Proof of
Theorem 1.1}. Let $s\in\{1,2\}$ and let $x$ be any $p$-adic integer. We claim that
$$\da_{s,2}\,\da_{p,3}+\sum_{n=1}^{p-1}\f{D_n(x)}{n^s}\eq\sum_{k=1}^{(p-1)/2}\f{(-x)^k}{k^s}\ (\mo\ p).\tag3.4$$
Clearly,
$$\sum_{n=1}^{p-1}\f{D_n(x)-1}{n^s}=\sum_{n=1}^{p-1}\f{\sum_{k=1}^n\bi{n+k}{2k}\bi{2k}kx^k}{n^s}
=\sum_{k=1}^{p-1}\bi{2k}kx^k\sum_{n=k}^{p-1}\f{\bi{n+k}{2k}}{n^s}.$$
Note that $\sum_{n=1}^{p-1}1/n^s\eq -\da_{s,2}\,\da_{p,3}\ (\mo\ p)$ since
$$\sum_{k=1}^{p-1}\f1{(2k)^s}\eq\sum_{n=1}^{p-1}\f1{n^s}\ (\mo\ p).$$
As $p\mid\bi{2k}k$ for $k=(p+1)/2,\ldots, p-1$, and
$$\sum_{n=k}^{p-1}\f{\bi{n+k}{2k}}{n^s}=\sum_{r=0}^{p-1-k}\f{\bi{2k+r}r}{(k+r)^s}
\eq(-2)^s\sum_{r=0}^{p-1-k}\f{(-1)^r\bi{p-1-2k}r}{(p-2k-2r)^s}\ (\mo\ p)$$
for $k=1,\ldots,(p-1)/2$, by applying Lemma 3.1 we obtain from the above that
$$\align \da_{s,2}\,\da_{p,3}+\sum_{n=1}^{p-1}\f{D_n(x)}{n^s}
\eq&(-2)^s\sum_{k=1}^{(p-1)/2}\bi{2k}kx^k\f{(-16)^{(p-1)/2-k}}{(p-2k)^s\bi{p-1-2k}{(p-1)/2-k}}
\\\eq&\sum_{k=1}^{(p-1)/2}\bi{2k}k\f{x^k}{k^s}4^{(p-1)/2-k}\bi{-1/2}{(p-1)/2-k}^{-1}
\\\eq&\sum_{k=1}^{(p-1)/2}\bi{2k}k\f{x^k}{k^s4^k}\bi{(p-1)/2}k^{-1}
\\\eq&\sum_{k=1}^{(p-1)/2}\bi{2k}k\f{x^k}{k^s4^k}\bi{-1/2}k^{-1}=\sum_{k=1}^{(p-1)/2}\f{(-x)^k}{k^s}
\ (\mo\ p).\endalign$$

In the case  $s=2$ and $x=1$, (3.4) yields the congruence
$$\da_{p,3}+\sum_{n=1}^{p-1}\f{D_n}{n^2}\eq\sum_{k=1}^{(p-1)/2}\f{(-1)^k}{k^2}\ (\mo\ p).$$
By Lehmer [L, (20)],
$$\sum^{(p-1)/2}\Sb k=1\\2\mid k\endSb\f1{k^2}\eq\da_{p,3}+\l(\f{-1}p\r)E_{p-3}\ (\mo\ p)$$
and hence
$$\sum_{k=1}^{(p-1)/2}\f{(-1)^k}{k^2}=2\sum^{(p-1)/2}\Sb k=1\\2\mid k\endSb\f1{k^2}-\sum_{k=1}^{(p-1)/2}\f1{k^2}
\eq\da_{p,3}+2\l(\f{-1}p\r)E_{p-3}\ (\mo\ p)$$
since $\sum_{k=1}^{(p-1)/2}(1/k^2+1/(p-k)^2)=\sum_{k=1}^{p-1}1/k^2\eq0\ (\mo\ p)$ if $p>3$.
So (1.1) follows.

With the help of (3.4) in the case $s=x=1$, we have
$$\align\sum_{n=1}^{p-1}\f{D_n}n\eq&\sum_{k=1}^{(p-1)/2}\f{(-1)^k}k
\eq\f12\sum_{k=1}^{(p-1)/2}\l(\f{(-1)^k}k+\f{(-1)^{p-k}}{p-k}\r)
\\\eq&-\f12\sum_{k=1}^{p-1}\f1k\bi{p-1}{k-1}=-\f1{2p}\sum_{k=1}^{p-1}\bi pk=-q_p(2)\ (\mo\ p).
\endalign$$
This proves (1.2).

Now fix a $p$-adic integer $x$.
Observe that
$$\align p\sum_{k=1}^{(p-1)/2}\f{(-x)^k}k\eq&-2\sum_{k=1}^{(p-1)/2}\f p{2k}\bi{p-1}{2k-1}(-x)^k
\\=&\sum^p\Sb j=1\\2\mid j\endSb \bi pj(-1)^{p-j}((\sqrt{-x})^j+(-\sqrt{-x})^j)
\\=&(-1+\sqrt{-x})^p+(-1-\sqrt{-x})^p+2\ (\mo\ p^2).
\endalign$$
Combining this with (3.4) in the case $s=1$ we immediately get (1.3).

The proof of Theorem 1.1 is now complete. \qed

\Remark\ 3.1. By modifying our proof of (1.2) and using the new identity
$\sum_{r=0}^{2n}\bi{2n}r/(2n+1-2r)=2^{2n}/(2n+1)$, we can prove the congruence
$\sum_{k=1}^{p-1}(-1)^ks_k/k\eq 4((\f 2p)-1)\ (\mo\ p)$ for any odd prime $p$.
Combining this with $\sum_{k=1}^{p-1}(-1)^kD_k/k\eq -4P_{p-(\f2p)}/p\pmod p$
(an equivalent form of (1.7)) we obtain that
$\sum_{k=1}^{p-1}(-1)^kS_k/k\eq 4(1-(\f 2p)-P_{p-(\f2p)}/p)\pmod p.$

\medskip
\noindent{\it Proof of Corollary 1.1}. Note that $\omega=(-1+\sqrt{-3})/2$ is a primitiv
e cubic root of unity.
If $p\not=3$, then
$$(-1+\sqrt{-3})^p+(-1-\sqrt{-3})^p=(2\omega)^p+(2\omega^2)^p=-2^p$$
and hence (1.3) with $x=3$ yields the congruence in (1.4).

Clearly (1.5) follows from (1.3) with $x=-4$.

Since $2^p-4^p+2=(2-2^p)(2^p+1)\eq 6(1-2^{p-1})\ (\mo\ p^2)$, (1.3) in the case $x=-9$ yields (1.6).

The companion sequence $\{Q_n\}_{n\gs0}$ of the Pell sequence is defined by $Q_0=Q_1=2$
and $Q_{n+1}=2Q_n+Q_{n-1}\ (n=1,2,3,\ldots)$. It is well known that
$$Q_n=(1+\sqrt2)^n+(1-\sqrt2)^n\quad\ \t{for all}\ n\in\N.$$
(1.3) with $x=-2$ yields the congruence
$$\sum_{k=1}^{p-1}\f{D_k(-2)}k\eq\f{2-Q_p}p\ (\mo\ p).$$
Since $Q_p-2\eq 4P_{p-(\f 2p)}\ (\mo\ p^2)$ by the proof of [ST, Corollary 1.3],
(1.7) follows immediately.

Recall that the Lucas sequence $\{L_n\}_{n\gs0}$ is given by
$$L_0=2,\ L_1=1,\ \t{and}\ L_{n+1}=L_n+L_{n-1}\ (n=1,2,3,\ldots).$$
It is well known that
$$L_n=\l(\f{1+\sqrt5}2\r)^n+\l(\f{1-\sqrt5}2\r)^n\ \t{for all}\ n\in\N.$$
Putting $x=-5$ in (1.3) we get
$$\align \sum_{k=1}^{p-1}\f{D_k(-5)}k\eq&\f{2-2^pL_p}p=\f{2^p(1-L_p)+2-2^p}p
\\\eq&-\f2p(L_p-1)-2q_p(2)\ (\mo\ p).
\endalign$$
It is known that $2(L_p-1)\eq 5F_{p-(\f p5)}\ (\mo\ p^2)$ provided $p\not=5$ (see the proof of [ST, Corollary 1.3]).
So (1.8) holds if $p\not=5$. We are done. \qed

\enddocument